\documentclass[12pt,reqno]{amsart}
\usepackage {amssymb}
\usepackage {amsmath}
\usepackage {bbm}
\usepackage{amsthm}
\usepackage{mathtools}
\usepackage{graphicx}
\usepackage {amscd}
\usepackage[alphabetic]{amsrefs}
\usepackage[colorlinks, citecolor=blue]{hyperref}
\usepackage{enumerate}
\usepackage{tikz}
\usepackage{tikz-cd}
\usepackage{verbatim,color,geometry}
\usepackage[all]{xy}
\usepackage{enumitem}
\usepackage{chngcntr}
\newtheorem{prop}{Proposition}
\newtheorem{theorem}[prop]{Theorem}

\newtheorem{problem}[prop]{Problem}
\newtheorem{conjecture}[prop]{Conjecture}

\theoremstyle{definition}

\newtheorem{rem}[prop]{\it Remark}

\newtheorem{say}[prop]{}
\newtheorem{lem-defn}[prop]{Lemma-Definition}
\newtheorem{defn-lem}[prop]{Definition-Lemma}

\newtheorem*{claim*}{Claim}

\newcommand{\bG}{\mathbb{G}}

\newcommand{\hvol}{\widehat{\rm vol}}

\newcommand{\cF}{\mathcal{F}}

\newcommand{\cE}{\mathcal{E}}

\newcommand{\fa}{\mathfrak{a}}

\newcommand{\fm}{\mathfrak{m}}

\newcommand{\Spec}{\mathrm{Spec}~}

\newcommand{\mult}{\mathrm{mult}}

\newcommand{\ord}{\mathrm{ord}}

\newcommand{\Gr}{\mathrm{Gr}}

\newcommand{\Val}{\mathrm{Val}}

\newcommand{\rd}{\mathrm{d}}

\newcommand{\ZZ}[1]{{\textcolor{blue}{[ZZ: #1]}}}

\numberwithin{equation}{section}

\begin{document}
\counterwithout{equation}{section}

\title{Open problems in K-stability of Fano varieties}

\date{}

\author{Chenyang Xu}
\address{Department of Mathematics, Princeton University, Princeton, NJ 08544, USA}
\email     {chenyang@princeton.edu}

\author{Ziquan Zhuang}
\address{Department of Mathematics, Johns Hopkins University, Baltimore, MD 21218, USA}
\email{zzhuang@jhu.edu}

\maketitle

\begin{abstract}
 In this note, we discuss a number of open problems in K-stability theory.   
\end{abstract}

\setcounter{tocdepth}{1}
\section{Introduction}
\subsection{K\"ahler-Einstein problem}
It has been a central question in geometry to find an `optimal' structure on a manifold $X$. In differential geometry, metrics $g$ satisfying the vacuum Einstein field equations, i.e. ${\rm Ric}(g)=\lambda g$, provide a candidate, and have been investigated extensively. When $X$ is a complex manifold, then it is also natural to consider focusing on  K\"ahler metrics $g$. For any K\"ahler metric $g$, if we consider the associated K\"ahler form $\omega$, then the Einstein equation becomes 
\begin{equation}\label{eq-KE}
{\rm Ric}(\omega)=\lambda \omega
\end{equation} 
and we can assume $\lambda=-1,0$ or 1. If we start with a background K\"ahler metric $\omega_0$ such that $\lambda \omega_0\in c_1(X)$, using the $\partial\overline{\partial}$-lemma, we can write ${\rm Ric}(\omega_0)-\lambda \omega_0=i \partial\overline{\partial} F$, and $\omega-\omega_0=i\partial\overline{\partial} \varphi$. Then the K\"ahler-Einstein equation \eqref{eq-KE} becomes a complex Monge-Amp\`ere equation
\begin{equation}\label{eq-MA}
{\displaystyle (\omega _{0}+i\partial {\bar {\partial }}\varphi )^{n}=e^{F-\lambda \varphi }\omega _{0}^{n}} \,. 
\end{equation}
Its behaviour is highly sensitive to the sign of the topological constant $\lambda$. It has been proved that when $\lambda=-1$ or 0, there is always a solution. However, when $\lambda=1$, the answer is much more delicate. We note in this case, $c_1(X)$ contains a K\"ahler class, i.e. the canonical class $-K_X$ is ample, and we call such $X$ Fano manifolds.

To solve $\varphi$, we look at the space of K\"ahler potentials
\[
\mathcal{H}=\{\varphi \mid \omega_0+i\partial\overline{\partial} \varphi >0 \} \, ,
\]
and it is shown in \cite{Mabuchi-Kenergy} there exists a K-energy functional $\mathcal{K}\colon \mathcal{H}\to \mathbb R$, such that a solution of \eqref{eq-MA} corresponds to a critical point of $\mathcal{K}$. Then in \cite{Tian-97} (also see \cite{DT-92}), a key speculation is made that to test whether there is a critical point, it should suffice to look at the slope of $\mathcal{K}$ along the ray of metrics arising from the pull back of the (normalized) Fubini-Study metic along embeddings $\lambda\cdot X\subset \mathbb{P}^N := \mathbb{P}(H^0(rK_X)^*)$ for all one parameter subgroups $\lambda \colon \mathbb{C}^*\to {\rm PGL}(N+1)$ and $r\gg 0$. Such a condition is coined as {\it K-stability}. Later, in \cite{Don-K}, it is realized that the K-stability condition is indeed algebraic and can be generalized to any polarized projective schemes.   With many researchers' contribution, now it is known that  K-polystability of a Fano variety is equivalent to the existence of a K\"ahler-Einstein metric, first in the smooth case (see e.g. \cites{CDS,Tian-KE}), then in the general case (see e.g. \cite{BBJ-uk=ke, Li-uk=ke, LXZ-HRFG}).

\subsection{K-stability of Fano varieties in algebraic geometry}There has been tremendous progress in the study of the K\"ahler-Einstein problem, from many different perspectives. One remarkable connection that happened in the last decade is its deep interaction with higher dimensional geometry.  As a result, algebraic K-stability theory of Fano varieties itself has grown to be a major field in algebraic geometry. 

Among all achievements of K-stability theory, the most striking one is the construction of projective scheme  parametrizing K-polystable Fano varieties with fixed numerical invariants, see Theorem \ref{K-moduli}. This construction involves many recipes: different ways to characterize K-stability (e.g. Theorem \ref{thm-FujitaLi}), progress on minimal model program (e.g. \cite{Bir-BABI}, Theorem \ref{thm-HRFG}), local KSBA stability of families of higher dimensional varieties (\cite{Kollar-modulibook}), and abstract stack theory (\cite{AHH-goodmoduli}) etc..   

Another surprising aspect is that the local-to-global correspondence between Fano varieties and klt singularities can be extended to include stability theory. More precisely, there is a local stability theory, developed based on the notion of normalized volume functions (see Paragraph \ref{say-normalizedvolume} and Theorem \ref{thm-SDC}). Another deep feature is the boundedness, see Theorem \ref{thm-localbounded}.

There have been several papers surveying the results in K-stability, see e.g. \cites{LLX-localsurvey, Xu-survey,Zhuang-localsurvey}. For readers who want to know the topic in more details, we refer to \cite{Xu-book}. In this note, we would like to discuss a number of open questions, rising from different aspects of the theory. After the construction of K-moduli spaces and the proof of the stable degeneration theorem, we hope these questions can lead to future development of K-stability theory, especially by establishing more connections with other fields of higher dimensional geometry. 

\vspace{2mm}

\noindent {\bf Acknowledgement}: We would like to thank Baohua Fu, Adrian Langer and Jie Liu for helpful discussions, and Christopher Hacon as well as the anonymous referee for useful suggestions. We are especially grateful to Kento Fujita to send us the calculation in \cite{Fujita26} for \ref{ex-contactexample}. CX is partially supported by NSF Grant DMS-2201349 and a Simons Investigator grant. ZZ is partially supported by the NSF Grants DMS-2240926, DMS-2234736, a Sloan research fellowship and a Packard fellowship. CX and ZZ are also partially supported by the Simons Collaboration Grant on Moduli of Varieties.

\section{K-stability of projective varieties}

\begin{say}[Stability threshold] One key feature of K-stability for Fano varieties is that there are many equivalent ways to define them, starting from the original definitions by looking at Futaki invariants ${\rm Fut}(\mathcal{X})$ of test configurations $\mathcal{X}$ (see \cites{Tian-97, Don-K}), and later by Ding invariants ${\rm Ding}(\mathcal{X})$ (see \cite{Ber-keimpliesk}), Fujita-Li valuative criterion \cites{Fuj-valuative, Li-minimizer}, and using Ding invariants ${\mathbf{D}}(\mathcal{F})$ of filtrations $\mathcal F$ \cite{Fuj-valuative, XZ-CMpositive, BLXZ-soliton}. 

From the computational view point, one effective way to verify K-stability of a log Fano pair $(X,\Delta)$ is using the \emph{stability threshold}
\begin{equation}\label{eq-delta}
\delta(X,\Delta)=\inf_E\frac{A_{X,\Delta}(E)}{S(E)} \, .
\end{equation}
Here $E\subset Y$ is a divisor on a normal birational model $\mu\colon Y\to X$, $A_{X,\Delta}(E)$ is the log discrepancy of $E$ with respect to $(X,\Delta)$,
\[
S(E)=\frac{1}{(-K_X-\Delta)^n}\int {\rm vol}(\mu^*(-K_X-\Delta)-tE) \, {\rm d}t \, ,
\]
 and the infimum runs through all divisors $E$ over $X$. The invariant $\delta(X)$ is defined in \cite{FO-delta}, and the equivalent description \eqref{eq-delta} is established in \cite{BJ-delta}.
\end{say}
\begin{theorem}\label{thm-FujitaLi}
We have the following: A log Fano pair $(X,\Delta)$ is 
\begin{enumerate}
\item (\cites{Fuj-valuative, Li-minimizer}) K-semistable if and only if $\delta(X,\Delta)\ge 1$; and
\item (\cite{LXZ-HRFG}) K-stable if and only if $\delta(X,\Delta)>1$.
\end{enumerate}
\end{theorem}
\begin{say}
To distinguish strict K-semistablity and K-polystability, we can follow the same idea and define reduced $\delta$-invariant $\delta^{\rm red}$ of a semistable log Fano pair $(X,\Delta)$ with a torus $\mathbb T$ action, as in \cite{XZ-CMpositive}. Then K-polystability of $(X,\Delta)$, which is equivalent to that any divisorial valuations $E$ satisfying $A_{X,D}(E)=S(E)$ arises from some $\bG_m$ action, can be characterized by  $\delta^{\rm red}(X,\Delta)>1$ by \cite{LXZ-HRFG}.
\end{say}

\begin{theorem}[K-moduli]\label{K-moduli}
Fix two integers $n$, $N$ and a positive number $V$. 
\begin{enumerate}
\item The functor 
\[
\mathfrak{X}_{n,N,V}^{\rm K}(S)
=
\left\{
\begin{array}{l}
\text{$(X,\Delta)\to S$ such that $X\to S$ is projective flat with}\\[0.2em]
\text{normal fibers, $\omega^{[m]}_{X/S}(m\Delta)$ is flat over $S$ for any $m$}\\[0.2em] 
\text{divisible by $N$, $N\Delta$ is integral and K-flat \cites{Kollar-modulibook},}\\[0.2em]
\text{fibers $(X_s,\Delta_s)$ are K-semistable log Fano pairs,}\\[0.2em]
\text{$\dim(X_s)=n$, $(-K_{X_s}-\Delta_s)^n=V$}
\end{array}
\right\}
\]
is represented by a finite type Artin stack $\mathfrak{X}^{\rm K}_{n,V,N}$.
\item $\mathfrak{X}_{n,N,V}^{\rm K}$ admits a separated good moduli space $\mathfrak{X}_{n,N,V}^{\rm K}\to {X}_{n,N,V}^{\rm K}$.
\item ${X}_{n,N,V}^{\rm K}$ is proper.
\item There exists a CM line bundle $\lambda_{\rm CM}$ which is ample on ${X}_{n,N,V}^{\rm K}$.
\end{enumerate}
\end{theorem}
In (1), the boundedness was first proved in \cite{Jia-bounded} and a different proof was later given in \cite{XZ-uniqueness}; the openness of K-semistability was proved in \cites{Xu-quasimonomial, BLX-openness}. Using the criterion established in \cite{AHH-goodmoduli}, the existence of the separated good moduli space (2) was proved in \cite{BX-uniqueness, ABHX-reductivity}. The proof of (3) is based on a higher rank finite generation theorem proved in \cite{LXZ-HRFG}, then the properness is first confirmed  in \cite{BHLLZX} by constructing a $\Theta$-stratification on the moduli of all Fano varieties $\mathfrak{X}_{n,N,V}^{\rm Fano}$ and then apply the general semistable reduction theorem in $\Theta$-stratification theory established in \cite{HL-instability} (see also \cite{AHH-goodmoduli}). Later a proof replacing the use of the $\Theta$-stratification by purely birational geometry arguments is given in \cite{BLXZ-properness}.  Finally, (4) is proved in \cites{CP-CMpositive, XZ-CMpositive}.

\begin{say}
It is proved in \cite{Zhuang-equivariant} the notions of K-stability and K-polystability do not depend on the ground field, so $\mathfrak{X}^{\rm K}_{n,V,N}$ and ${X}_{n,N,V}^{\rm K}$  are well defined moduli spaces over $\mathbb Q$.

One may naturally wonder whether one can extend the K-moduli construction to fields of positive characteristics, or even in the arithmetic setting, e.g. over $\mathbb {Z}$. For now some fundamental challenge remains. Most importantly, a large part of birational geometry known in characteristic 0, e.g. resolution of singularities and MMP theory, remains wide open in the more general case, except low dimensional cases. 

However, a small part of K-stability theory in general dimension is extended to the arithmetic case \cites{Berman-height,AB-heightI,AB-heightII}. This presents evidence that there should be an interesting theory also in this context.
\end{say}

\subsection{K-stability of explicit Fano varieties}

 Historically, the question of verifying whether a given Fano admits a K\"ahler-Einstein metric has been well-studied. The algebraic version, which  currently presents the best approach to the original question, is to investigate K-(semi,poly)stability of any given Fano variety (or more generally log Fano pair). This question has drawn a lot of interest, partly because it connects to many existing techniques developed in extensive research on Fano varieties. Here we discuss several classes of examples of particular interest to us (also see Problem \ref{prob-slodowy} for other examples). 

\subsubsection{Fano hypersurfaces}
\begin{problem}\label{prob-hypersurface}
For any $3\le d \le n+1$, prove every degree $d$ smooth hypersurface in $\mathbb P^{n+1}$ is K-stable. 
\end{problem}

When $d=n+1$ this is first proved in \cite{Fujita-index1}. Later the Abban-Zhuang method is invented in \cites{AZ-1,AZ-2}, which in particular covers the case 
\[
n+1-\max\{1,n^{\frac 13}\}\le d\le n+1\, .
\]It is also easy to see that a general degree $d$ $(d\ge 3)$ hypersurface is K-stable by looking at the Fermat hypersurface. So there is a K-moduli space $H^{\rm K}_{n,d}$ whose general points parametrize general $n$-dimensional hypersurfaces of degree $d$ $(3\le d\le n+1)$. Problem \ref{prob-hypersurface} has an affirmative answer for $n\le 4$ by combining these works with the case $d=3$ proved in \cite{LX-cubic, Liu-cubic}.

\begin{problem}
Describe the K-moduli space $H^{\rm K}_{n,d}$ $(3\le d\le n+1)$.  
\end{problem}
This problem seems out of reach for now, in general. Currently, only when $d=3$, we have a conjectural picture that the K-moduli $H^{\rm K}_{n,3}$ is the same as the GIT moduli $H^{\rm GIT}_{n,d}$ of cubic hypersurfaces, and this is confirmed for $n\le 4$. For $d\ge 4$, the picture is very unclear (see \cite{ACKLP-quartic} for some partial description of $H^{\rm K}_{3,4}$).

It is also interesting to look at higher codimension Fano complete intersections. 
\begin{problem}
For any $d_1,\dots, d_k$ with $k+2\le \sum_i d_i\le n+k$, show that a general complete intersection of $(d_1,\dots,d_k)$ is K-stable.
\end{problem}

We can also look at hypersurfaces in weighted projective spaces, see \cite{ST-weighthypersurface} for partial results. 

\subsubsection{Moduli of vector bundles}

Another natural Fano variety is the moduli space $M_{d,r}(C)$ of rank $r$ polystable vector bundles on a fixed smooth projective curve $C$, with a fixed determinant $L$ of degree $d$. $M_{d,r}(C)$ is a Fano variety and smooth when ${\rm gcd}(r, d)=1$.

\begin{problem}
\begin{enumerate} 
\item Show $M_{d,r}(C)$ is K-stable for any smooth projective curve $C$ of genus $g$.
\item The K-moduli space $X_g$ containing $[M_{d,r}(C)]$ admits a birational map to $\overline{M}_g$, which is conjecturally isomorphic over $M_g$ by (1). Give a modular interpretation of the limiting points on $X$.
\end{enumerate}
\end{problem}
Not much is known except in \cite{Keller-rank2}, it is proved that when $r=2$, $d=1$, $M_{1,2}(C)$ is K-stable for a curve $C$ corresponding to a general point in $M_g$.


\subsubsection{Low dimensional cases}
There has been a huge amount of literature devoted to study K-stability of low dimensional varieties. The case of $\dim(X)=2$ is essentially known (see e.g. \cite{Tian-surface,MM-degree4,OSS-surface}).  When $\dim(X)=3$, all smooth Fano threefolds are classified by Iskovskikh (when $\rho(X)=1$) and Mori-Mukai (when $\rho(X)\ge 2$), and there are 105 deformation families. Many families have been studied in depth, see e.g. \cite{LX-cubic,AZ-1, Fujita-3.11, threefoldsII, LiuZhao-threefold, CDF-picard2, 2.11} and many others, also see the webpage \cite{Fanography}. In the prime case, we want to ask the following question. 
\begin{problem}
Show that any smooth Fano threefold $X$ with $\rho(X)=1$ is K-semistable. 
\end{problem}
This is known by \cite{AZ-1}, except two families, namely, No $1.9$ and No $1.10$. The family No $1.10$ is especially subtle, since it contains smooth strictly K-semistable Fano manifolds, which degenerate to the Mukai-Umemera smooth threefold known to be K-polystable with a positive dimensional automorphism group (see \cite{Tian-97}). 

Another type of examples are Gushel-Mukai Fano manifolds of dimension $n$ ($3\le n\le 6$).

\begin{problem}
Show that any smooth Gushel-Mukai Fano variety $X$ is K-stable.  
\end{problem}
When $n=3$, this is answered in \cite{AZ-1}.  For $n=4,5,6$, this is only known for general members by the work of \cite{LW-GM}.

\subsection{Stability of tangent bundle}
It is known that K\"ahler-Einstein manifolds (for all possible values of   $\lambda$ in \eqref{eq-KE}) have a stable tangent bundle, because the K\"ahler-Einstein metric $g_{\rm KE}$ on the tangent bundle $T_X$ will make $(T_X,g_{\rm KE})$ a Hermitian-Einstein metric with respect to the K\"ahler form $\omega_{\rm KE}$. Then we can apply (the easier direction of) the Donaldson-Uhlenbeck-Yau Theorem which says a vector bundle $E$ admits a Hermitian-Einstein metric $g$ with respect to  a given K\"ahler form $\omega$ implies that $E$ is polystable. 

A folklore question, which will be an important step to understand stability, is to give a purely algebraic proof of these statements. 
\begin{problem} 
Let $X$ be a  Fano variety. Give an algebraic proof of the following:
\begin{enumerate}
\item If $X$ is K-semistable, then the tangent sheaf $T_X$ is semistable  with respect to the polarization $-K_X$. 
\item If $T_X$ is polystable, then there is a quasi-\'etale cover $f\colon Y\to X$ such that $Y=\prod_i X_i$ and $T_Y=f^*(T_X)=\bigoplus_i T_{X_i}$ such that $T_{X_i}$ is stable with respect to $-K_Y$.
\end{enumerate}
\end{problem}

A maximal destabilizing sheaf $F\subset T_X$ will induce a zero morphism
\[
[F,F] \to T_X/F \, ,
\]
which means it is an integrable foliation. Moreover, for a general complete intersection curve $C$ of $H_i\in |-mK_X|$ $(m\ge 0)$,  $F|_{C}$ is ample as it is a stable vector bundle of positive slope by the restriction theorem. Therefore, the foliation is indeed algebraic \cites{BM-foliation, KST-foliation}, i.e. it is given by  a rational map $X\dasharrow Z$.

For a similar statement  in the Calabi-Yau case ($\lambda=0$), (1) follows from the fact that the maximal destabilizing subsheaf $F$ yields an algebraically integrable foliation by the above argument, and moreover its leaves are rationally connected, a contradiction; nevertheless, the known proof of (2) is of an analytic nature, see \cite{HP-decomposition}. The canonically polarized case is completely open: as far as we know, there is no purely algebraic proof of (1) even for minimal general type surfaces.

In fact, following \cite{Tian-tangentstability} (also see \cite{DGP-stabletangent}), we can expect a stronger result by looking at the following construction: the image of the class $c_1(X)\in {\rm Pic}(X)\otimes \mathbb{Q}^*$ in $H^1(X,\Omega_X^{[1]})$ determines a canonical extension
\[
0\to  \Omega_X^{[1]}  \to \mathcal{E}\to \mathcal{O}_X\to 0 \, .
\]
\begin{problem}\label{problem-canonicalexten} 
Let $X$ be a  Fano variety. If $X$ is K-semistable (resp. K-polystable), give an algebraic proof of the following: the canonical extension $\mathcal{E}$ is semistable (resp. polystable) with respect to the polarization $-K_X$. 
\end{problem}

\subsection{Gibbs stability}

\begin{say}
In \cite{Berman-Gibbs}, an invariant $\gamma(X)$ is introduced to characterize a notion called \emph{Gibbs stability}, which is conjectured to be identical to $\delta(X)$ in the log Fano case. 

More precisely, for a Fano variety $X$, let $D_m$ be the divisor on $X^{N_m}:=X\times \cdots \times X$ ($N_m$ times), where $N_m=H^0(-mK_X)$,   given by the pulling back of the vanishing locus of ${\rm Det}(x_{ij})_{1\le i,j\le N_m}$ under the embedding $f^{N_m}\colon X^{N_m}\to (\mathbb{P}^{N_m-1})^{N_m}$. Then one can define an invariant
\begin{equation}\label{eq-gamma}
\gamma(X)=\liminf_m{\rm lct}(X^{N_m},\frac{1}{m}D_m)\, .
\end{equation}
This arises from a conjectural picture of constructing the K\"ahler-Einstein measure from a probability theoretic viewpoint: The condition $\gamma(X)>1$, i.e. \emph{Gibbs stable}, guarantees the measure 
\[
\nu^{(N_m)}:= i^{N_m^2} ({\rm det}(S^{(N_m)}))^{-\frac{1}{m}} \wedge \overline{({\rm det}(S^{(N_m)}))^{-\frac{1}{m}}}
\] on $X^{N_m}$ is of finite mass for $N_m\gg 1$, where $S^{(N_m)}$ is the section in $H^0(-mK_{X^{N_m}})$ corresponding to $D_m$. So under this condition, one can define the probability measure on $X^{N_m}$
\[ 
\mu^{(N_m)}=\frac{1}{Z_{N_m}}\nu^{(N_m)},  \ \ \mbox{where } Z_{N_m}:=\int_{X^{N_m}}\nu^{(N_m)}<+\infty\, ,
\]
which yields a random point process on $X$ with $N_k$ particles.

The empirical measure of a random point process is
\[
 \delta_{N_m}:=\frac{1}{N_m}\sum^{N_m}_{i=1} \delta_{x_i} \,,
\]
which defines a map $\delta_{N_m}\colon X^{N_m}\to \mathcal{M}_1(X)$ where $\mathcal{M}_1(X)$ is the space of normalized measures of $X$.
Then the laws of the empirical measures for the measure $\mu^{(N_m)}$ on $X^{N_m}$ (conjecturally) satisfies a large deviation principle with speed $N_m$ and rate functional $F(\mu)$. Here $F(\omega^n/V)$ is Mabuchi’s K-energy of the K\"ahler form $\omega$ (normalized so that $F$ vanishes on the normalized measure ${\rm d}V_{\rm KE}$ induced by the K\"ahler-Einstein metric). Roughly speaking, this is 
\[ 
{\rm Prob}_{X^{N_m},\mu^{(N_m)}}\left(\frac{1}{N_m}\sum_i \delta_{x_i}\in B_{\varepsilon}(\mu) \right)\sim e^{-N_mF (\mu)}
\] 
as first $N_m \to +\infty$ and then $\varepsilon \to 0$, where $B_{\varepsilon}(\mu)$ denotes a ball of radius $\varepsilon$ centered at a given probability measure $\mu$ in the space $\mathcal{M}_1(X)$ of all normalized measures $\mu$ on $X$. 

As a consequence, the empirical measures of the canonical random point processes (conjecturally) converges in probability to the normalized volume form ${\rm d}V_{\rm KE}$ for the K\"ahler-Einstein metric $\omega_{\rm KE}$ (see \cites{Berman-Gibbs, ABS-Gibbs}).  

A similar picture of this perspective is established for $X$ with ample $K_X$. For Fano manifolds it is largely conjectural, for instance the following is not known.
\begin{conjecture}[Berman]
For a Fano variety $X$, $\gamma(X)>1$ (resp. $\ge 1$) if and only if $X$ is K-stable (resp. K-semistable).
\end{conjecture}
One also expects $\delta(X)=\gamma(X)$ when $\delta(X)\le 1$. One direction $\delta(X)\ge \gamma(X)$ is known (see \cite{Fujita-Gibbs}, \cite{FO-delta}). 
\end{say}


\section{Local K-stability}
\begin{say}\label{say-normalizedvolume}
The starting point of the local K-stability theory is the definition of the normalized volume of valuations as in \cite{Li-minimizer}. Let $x\in (X={\rm Spec}(R),\Delta)$ be a klt singularity germ of dimension $n$. Let ${\rm Val}_{X,x}$ be the space of all valuations whose center is $x\in X$. For any $v\in {\rm Val}_{X,x}$, we can define its log discrepancy $A_{X,\Delta}\colon {\rm Val}_{X,x} \to (0, +\infty]$ and volume ${\rm vol}(v)=\lim_{k\to \infty}\frac{n!}{k^n}\ell (R/\mathfrak{a}_k)$ for $\mathfrak{a}_k=\{f\in R \mid v(f)\ge k \}$. Then the key definition is
\[
\widehat{\rm vol}(v) = 
\begin{cases}
A_{X,\Delta}(v)^n\cdot {\rm vol}(v), & \text{if } A_{X,\Delta}(v) < +\infty,\\
+\infty,  & \text{if } A_{X,\Delta}(v) = +\infty;
\end{cases}
\]
and then we define \emph{volume} of a klt singularity $x\in (X,\Delta)$ to be 
\begin{equation}\label{df-volume}
\widehat{\rm vol}(x, X,\Delta)=\inf_{v\in  {\rm Val}_{X,x} }\widehat{\rm vol}(v) \,.
\end{equation}
\end{say}

A deep and surprising phenomenon is that \emph{any} klt singularity admits a canonical K-semistable degeneration, constructed using the minimizer of $\widehat{\rm vol}(v) $.

\begin{theorem}[Stable degeneration]\label{thm-SDC}  For any klt singularity $x\in (X,\Delta)$, we have the following
\begin{enumerate}
\item There exists a quasi-monomial minimizer $v$ of $\widehat{\rm vol}(v) $, and the minimizer is unique up to rescaling. 
\item The graded ring ${\rm Gr}_vR$ is finitely generated. 
\item Let $X_0={\rm Spec}({\rm Gr}_vR)$, $\Delta_0$ the divisorial degeneration of $\Delta$, $\xi_v$ the corresponding Reeb vector, then $(X,\Delta,\xi_v)$ is a K-semistable log Fano cone. 
\end{enumerate}
\end{theorem}
For (1), the existence of a quasi-monomial valuation is proved in \cites{Blu-existence, Xu-quasimonomial}, and the uniqueness is proved in \cites{XZ-uniqueness, BLQ-convexity}. The most difficult part, namely the finite generation (2) is proved in \cite{XZ-HRFG} (see also \cite{Chen-HRFG}). The K-semistability is addressed in \cites{Li-K=M, LL-K=M, LX-SDCI,LX-SDCII}. 

\begin{rem}
Together with \cite{LWX-tangent}, this produces a two-step canonical degeneration of any klt singularity to a K-polystable Fano cone singularity. Such two-step degeneration process, noticed first in \cite{DS-degeneration2}, also exists in other context in K-stability theory. For instance, using $H$-functional arising from the study of K\"ahler-Ricci soliton, \cite{HL-soliton } and \cite{BLXZ-soliton} produces a two-step degeneration of any Fano variety $X$ to a K-polystable pair $(Y,\xi)$, where $Y$ is a Fano variety with a torus $\mathbb T$-action and $\xi \in {\rm Hom}(\mathbb{G}_m,\mathbb{T})\otimes {\mathbb R}$.
\end{rem}

\subsection{Symplectic singularities and contact Fano varieties}

A surprising application of the general framework of the Stable Degeneration Theorem \ref{thm-SDC} is to investigate symplectic singularities. This is initiated in the recent work \cite{NO-symp}, which relies on analytic results. Here we want to present a more general purely algebraic conjectural picture. 
\begin{say}[Symplectic Singularity \cite{Beauville-singularity}]
A normal singularity $x\in X$ is called \emph{symplectic}  if its smooth part $X^{\rm reg}=X\setminus {\rm Sing}(X)$ carries a closed non-degenerate $2$-form $\omega$ whose pullback to any resolution $X'\to X$ extends to a holomorphic $2$-form on $X'$. If the symplectic singularity $(X,\omega)$ admits a good $\mathbb{G}_m$-action (i.e. $x$ is contained in the closure of any $\mathbb{G}_m$-orbit), then there exists a positive integer $w$ such that $\lambda\cdot \omega=\lambda^{w}\omega$. We call $x\in X$ a \emph{conical} symplectic singularity, and $w$ the {\it weight} of the $\mathbb G_m$-symplectic singularity $(X,\omega)$.
\end{say}

When applying the Stable Degeneration Theorem \ref{thm-SDC} to symplectic singularities, the following striking phenomenon seems to hold.  

\begin{problem}\label{p-symp}
For a $2n$-dimensional symplectic  singularity $x\in X$, and any symplectic form $\omega_X$, prove the stable degeneration $x_0\in X_0$ is also symplectic with the symplectic form $\omega_{X_0}$ the degeneration of $\omega_X$.
\end{problem}

\begin{say}[Analytic work] The primary evidence of Problem \ref{p-symp} is from \cite{NO-symp}.
Let $X$ be a singular hyperk\"ahler variety. Assume $X$ admits a crepant resolution, which implies it has a smooth hyperk\"ahler deformation $X_t$. Let $g_t$ be the hyperk\"ahler metric on $X_t$. Then the metric $g_t$ converge, in terms of Gromov-Hausdorff limit, to a hyperk\"ahler metric $g$ on $X$. For any $x\in (X,\omega)$, one can construct the metric tangent cone $o\in C$ as in \cite{DS-degeneration2}, which is the analytic precursor of Theorem \ref{thm-SDC}.

Then using that the symplectic forms are flat with respect to the metrics, it is shown in \cite{NO-symp} that in this case $(X_0,\omega_0)$ is also symplectic, answering Problem \ref{p-symp} in this case. Moreover,  \cite{NO-symp} conjectures the existence of a `good' Ricci-flat metric on any klt singularity germ, which will imply Problem \ref{p-symp} always holds by the same argument. 


\end{say}

Using properties of symplectic deformations of symplectic singularities, Problem \ref{p-symp} implies that $x\in X$ and $x_0\in X_0$ are analytically isomorphic (see \cite{NO-symp}). This gives a more refined conjectural picture of Kaledin's conjecture that any symplectic singularity is analytically conical (see \cite{Kaledin-symsing}).

 Global objects closely related to symplectic singularities are contact Fano varieties. 

\begin{say}[Contact Fano varieties]
A \emph{contact} variety is an algebraic variety $X$  of odd dimension $2n-1$  with klt singularities and an ample $\mathbb {Q}$-Cartiter reflexive rank one sheaf $L$ (i.e. for some positive integer $m$, the reflexive power $L^{[m]}$ is Cartier), such that on the smooth locus $X^{\rm reg}$, there exists a section $\theta\in H^0(\Omega^1_{X^{\rm reg}}\otimes L)$ yielding an exact sequence of vector bundles:
\[
0 \to F \to T_{X^{\rm reg}} \xrightarrow[\text{}]{\theta}  L|_{X^{\rm reg}} \to 0 \, ,
\]
which defines contact structure on $X^{\rm reg}$, i.e. the induced map ${\rm d} 
\theta\colon \wedge^2 F \to  L|_{X^{\rm reg}}$ is nowhere
degenerate. Equivalently one can demand that $\theta \wedge ({\rm d}\theta)^{\wedge n-1}$ as an element of
$H^0(X^{\rm reg}, \omega_{X^{\rm reg}} \otimes L^{n}|_{X^{\rm reg}})$ has no zeroes. Since $K_X=-nL$ under these assumptions and $L$ is ample,  $X$ is a klt Fano variety. 

For our purpose, we also need the extension of this definition  to a klt Fano orbifold $\mathcal{X}$, see  \cite{Namikawa-contact}. 
\end{say}

When $X$ is a $(2n-1)$-dimensional contact Fano variety with $-K_X=nL$, the cone $(o\in Y)=C(X,H)$ for any $H$ with $L=-rH$ yields a $\mathbb{G}_m$-symplectic singularity of weight $r$ (see \cite{Smiech-contact}). Conversely, for a conical symplectic singularity, a direct calculation shows that the base is a klt contact Fano orbifold (see \cite{Namikawa-contact}).

One strong conjecture predicts all contact Fano manifolds are necessarily homogeneous spaces (see \cite{Lebrun-contact}, where the existence of a K\"ahler-Einstein metric is also explicitly conjectured). However, there are a lot more singular examples of contact Fano varieties. In particular, there exist easy examples of contact Fano varieties which are not K-semistable. For instance, a weighted projective space $\mathbb{P}(a_1,\cdots, a_{2n})$ with $a_1+a_2=\cdots=a_{2n-1}+a_{2n}$ is a contact Fano variety with respect to the form descent from $\rd z_1\wedge \rd z_2+\cdots +\rd z_{2n-1}\wedge \rd z_{2n}$, whereas it is K-semistable if and only if $a_1=\cdots=a_{2n}$. Nevertheless, we expect many `natural' contact Fano varieties are K-polystable. A class of `natural' contact Fano varieties is the base of conical singularities arising from Lie theory.

\begin{problem}\label{prob-slodowy}More precisely,
\begin{enumerate}
\item Is the projective base of a nilpotent orbit closure  K-polystable?
\item Is the projective base of a Slodowy slice  K-polystable?
\end{enumerate}
\end{problem}
In both cases, the weight of the `natural' conical action is 1 or 2.

\begin{say}[An example, \cite{Fujita26}]\label{ex-contactexample}

A family of non-Gorenstein explicit examples of projective bases of Slodowy slices are given in \cite{LNSvS}:
for each integer $k>1$, consider the degree $4k$ hypersurface in $\mathbb{P}(2,2,2,2k-1,2k-1)$ given by the explicit equation
\[
X_k:= \left((xz-y^2)^n+u^2x+2uvy+v^2z=0 \right)\subset \mathbb{P}(2,2,2,2k-1,2k-1) ,
\]
with $\deg(x)=\deg(y)=\deg(z)=2$ and $\deg(u)=\deg(v)=2k-1$. 

We give a sketch of the calculation in \cite{Fujita26} which shows that $X_k$ is K-polystable.  
\begin{claim*}
$X_k$ is K-polystable for any $k>1$.
\end{claim*}

The hypersurface $X_k$ admits a $G:={\rm SL}_2(\mathbb C)$-action:
$
A=\begin{pmatrix}
a & b\\
c &d 
\end{pmatrix}
 \mbox{ \ sends \ }  [x:y:z:u:v] 
$ to  $[x':y':z':u':v']$ where
\[
\begin{pmatrix}
x' & y'\\
y' &z'
\end{pmatrix}
=A^{T}\begin{pmatrix}
x & y\\
y &z
\end{pmatrix}A  \mbox{ \ and \ } \begin{pmatrix}
u' \\
v' 
\end{pmatrix}  = A^{-1} \begin{pmatrix}
u\\
v 
\end{pmatrix}   \, .
\]
In particular,  $X_k$ admits a $\mathbb{G}_m$-action by $\begin{pmatrix}
\lambda & 0\\
0 &\lambda^{-1}
\end{pmatrix}$:
\[
\lambda\cdot  [x:y:z:u:v]\to [\lambda^2 x:y :\lambda^{-2}z:\lambda^{-1}u\colon \lambda v] \, .
\]
Let $L=(x=y=z=0)\subset X_k$ and $C=\left((u=v=0)\cap X_k\right)$. Then $L$ and $C$ are closed $G$-orbits.
Let $Z\subset X_k$ be  a $G$-invariant irreducible closed subvariety. If $L\neq Z$, then $Z$ contains a point $p$ with $(x,y,z)_p\neq (0,0,0)$, which we may further assume $z\neq 0$. Since $\lambda\cdot  [x:y:z:u:v]\to [\lambda^4 x:\lambda^{2} y :z:\lambda^{2k-2}u\colon \lambda^{2k} v]$, so
$(0,0,1,0,0)\in \overline{G\cdot p}\subset Z$, which implies $C\subset Z$.

By \cite{Zhuang-equivariant}, we only need to check for any $G$-invariant divisorial valuations $D$, $\frac{A_{X_k}(D)}{S(D)}>1$. The closure of its center ${\rm Cent}_{X_k}(D)$ is either $L$ or contains $C$. So it suffices  to  get the estimate $\delta_L(X_k), \delta_C(X_k)>1$,  where for any $Z\subset X_k$, 
\[
\delta_Z(X_k)=\inf\left\{\frac{A_{X_k}(D)}{S(D)}\, \Big| \, Z\subset \mbox{ closure of }{\rm Cent}_{X_k}(D)\right\}\,.
\]

We first analyze divisors $D$ with $L={\rm Cent}_{X_k}(D)$. We work on the model $f_k\colon X_k'\to X_k$ which is the blow up along $L$ with an exceptional divisor $E_k$, and we will pick a flag on $X_k'$ to apply the Abban-Zhuang method (see \cite{AZ-1} or \cite[Section 4.5]{Xu-book}) to estimate $\delta_L(X_k)$. 

Denote by $\mathbb{P}=\mathbb{P}(2,2,2,2k-1,2k-1)$ and the blow up $f\colon \widetilde{\mathbb P}\to \mathbb{P}$ along $L$ with the exceptional divisor $E\cong \mathbb{P}_{xyz}^2\times \mathbb{P}^1_{uv}$. Then $\widetilde{\mathbb P}$ admits a bundle structure $g\colon\widetilde{\mathbb P} \to \mathbb{P}^2_{x,y,z}\cong \mathbb{P}^2$ with fiber isomorphic to $\mathbb{P}(2,2k-1,2k-1)\cong \mathbb{P}(1,1,2)$.  Let $H_1=g^*\mathcal{O}_{\mathbb P^2}(1)$ and $H_2=f^*\mathcal{O}(2k-1)$. We have 
\begin{equation}\label{eq-intersectionnumber}
H_1^3=0, \, H_1^2\cdot H_2^2=\frac{1}{2}, \, H_1\cdot H^3_2=\frac{(2k-1)}{4}, \, \mbox{and \ } H_1^4=\frac{(2k-1)^2}{8}\, .
\end{equation}
 Denote by $(a,b)$ the divisor class of $aH_1+bH_2$ on $\widetilde{\mathbb P}$, then $(a,b)|_E=(a,b)_{\mathbb{P}^2\times \mathbb{P}^1}$. Let $D_x$ be the class given by $x=0$ and $D_x'$ its birational transform. Then 
\[
f^*D_x=D_x'+\frac{1}{2k-1}E \mbox{ and }  [D_x']=(1,0), \ [f^*D_x]=(0,\frac{2}{2k-1}) \, .
\] 
So $[E]=(-(2k-1),2)$. Similarly we can compute $[K_{\widetilde{\mathbb P}}]=(2k-4,-4)$ and $[X_k']=(1,2)$. 
The equation of $E_k$ in $\mathbb{P}_{xyz}^2\times \mathbb{P}^1_{uv}$ is $(u^2x+2uvy+v^2z=0)$ which is isomorphic to $\mathbb{P}^1\times\mathbb{P}^1$. The pull back of $\mathcal{O}_{\mathbb{P}^2\times \mathbb{P}^1}(1,0)$ (resp. $\mathcal{O}_{\mathbb{P}^2\times \mathbb{P}^1}(0,1)$) to $E_k$ is $\mathcal{O}_{\mathbb{P}^1\times \mathbb{P}^1}(1,1)$ (resp. $\mathcal{O}_{\mathbb{P}^1\times \mathbb{P}^1}(0,1)$). So putting this together for $E_k=E\cap X'_k$, we conclude $A_{X_k}(E_k)=\frac{2}{2k-1}$ and $(K_{X'_k}+E_k)|_{E_k}=K_{E_k}$ (i.e. there is no different divisor).

 Let $\ell$ be a line on $\mathbb{P}^1\times\mathbb{P}^1$ with the class $\mathcal{O}(1,0)$ and $q\in \ell$. We can see
\begin{enumerate}
\item For $t_1\ge 0$, $L_{t_1}:=f^*\mathcal{O}(1)-t_1E_k$ is nef or pseudo-effective if and only if $t_1\in [0,\frac{1}{2(2k-1)}]$.
\item Let $L_{t_1,t_2}:=L_{t_1}|_{E_k}-t_2\ell$, then $[L_{t_1,t_2}]=(t_1(2k-1)-t_2,\frac{1}{2k-1}+t_1(2k-3))$. Fix $t_1\in [0,\frac{1}{2(2k-1)}]$, for $t_2\ge 0$,  $[L_{t_1}|_{E_k}-t_2\ell]$ is nef or pseudo-effective if and only if $t_2\in [0, t_1(2k-1)]$.
\item Let $L_{t_1,t_2,t_3}=L_{t_1,t_2}|_{\mathbb P^1}-t_3[q]$. Then $\deg_{\mathbb P^1}(L_{t_1,t_2,t_3})=\frac{1}{2k-1}+t_1(2k-3)-t_3$.
\end{enumerate}
For the multigraded linear series  arising from repeated refinement of $f_k^*\mathcal{O}(1)$ by the flag 
\[
(q\in \ell\subset E_k\subset X'_k)\,,
\] 
the above analysis implies their $S$-functions can be computed as integral of functions with coefficients coming from intersection numbers. 
As a result, we can apply the Abban-Zhuang inequality to obtain the local stability threshold estimate $\delta_L(X_k)> 1$.

We can similarly blow up $C\subset X_k$ to get $f'_k\colon Y_k'\to X_k$ with an exceptional divisor $F_k\cong \mathbb{P}^1\times \mathbb{P}^1$. For any $G$-equivariant divisorial valuation $D$ whose closed center on $X_k$ contains $C$, its closed center  on $Y_k'$  contains $C'$ dominating $C$. One can apply the Abban-Zhuang method for the flag $C'\subset F_k\subset Y'_k$, and conclude that $\frac{A_{X_k}(D)}{S(D)}>1$. So
$\delta_C(X_k)>1$.
\end{say}

Recall that $E$ is a Koll\'ar component over a klt singularity $x\in X$ if there exists a projective morphism $Y\to X$ such that $E\subset Y$ is mapped to $x$ with $(Y,E)$ plt and $-E$ ample over $X$. Then there exists a DM stack $\mathcal{Y}\to Y$ such that  $(\cE\subset \mathcal{Y})$ is the index one cover of the $\mathbb{Q}$-Cartier divisor $(E\subset Y)$.
The following problem investigates the minimizer of a symplectic singularity and provides a more refined version of  Problem \ref{p-symp}.

\begin{problem}\label{prob-contactandsym}
Let $x\in X$ be a symplectic singularity. 
\begin{enumerate}
\item (Rank one, \cite{NO-symp}) The minimizer of $\hvol(x,X)$ is given by a Koll\'ar component $E$.
\item (Weight) $A_X(E)=rn$ for $r=1$ or $2$. 
\item (Contact structure) Let $\pi\colon \mathcal{Y}\to Y\to X$  be the morphism extracting $E\subset Y$ as a Koll\'ar component with the index one cover $(\mathcal{E}\subset \mathcal{Y})$ of $(E\subset Y)$. Then we have $\pi^*(\omega_X) \in H^0(\mathcal{Y}, \Omega^2_{\mathcal{Y}}({\rm log} \mathcal{E})(-r\mathcal{E}))$. Moreover, its image $\theta=\partial(\pi^*(\omega_X))\in H^0(\mathcal{E},\Omega^1_{\mathcal{E}}\otimes \mathcal{O}_{\mathcal Y}(-r\mathcal{E})|_{\mathcal{E}})$ under the residue morphism 
\[
\partial \colon \Omega^2_{\mathcal{Y}}({\rm log} \mathcal{E})(-r\mathcal{E}) \to \Omega^1_{\mathcal{E}}(-r\mathcal{E}|_{\mathcal{E}})
\]
satisfies that $(\mathcal{E},\theta)$ is an orbifold contact structure with the contact line bundle $\mathcal{O}_{\mathcal Y}(-r\mathcal{E})|_{\mathcal{E}}$. 
\item (K-polystability) $\mathcal{E}$ is K-polystable, which means  the underlying log Fano pair $(E,\Delta_E)$ is K-polystable.
\end{enumerate}
In particular, we have the following:
\begin{enumerate}
  \setcounter{enumi}{4}
  \item Any symplectic singularity $x\in X$ admits a Koll\'ar component which is a contact Fano variety (in the orbifold sense).
\end{enumerate}
\end{problem}

In (1), the minimizer being a Koll\'ar component corresponds to the `quasi-regular' case in the literature on Sasaki-Einstein metrics. In fact, it is inspired by the speculation that the Sasaki-Einstein metric in this setting should be 3-Sasaki. The assertion in (2) that $r=1$ or $2$  means that we expect any conical symplectic singularity (analytically) admits a \emph{canonical} weight one or two action, which might not be the same as the original one. While this could be too optimistic, we do not know any counterexample to this. For ADE surface singularities $\mathbb{C}^2/G$ ($G\subset {\rm SL}_2(\mathbb C)$),  $r$ is given by $2/|I|$ where $I$ is the kernel of $G\to {\rm PGL}_2(\mathbb C)$.

\vspace{2mm}

Problem \ref{prob-contactandsym}, if true, says any conical symplectic singularity is a `hidden hyperk\"ahler cone'. Therefore, it  would imply the  package for hyperk\"ahler cone singularities should conjecturally exist on \emph{any} conical symplectic singularities. This would have many other consequences. In the rest of this section, we explicitly spell out some of the questions.

\begin{say}[Hamiltonian reduction] The canonical structure given by the minimizer interacts with many other aspects of symplectic singularities. Let $(X={\rm Spec}(R),\omega)$ be a conical symplectic singularity, and $\mathbb T$ a torus which admits an algebraic Hamiltonian action on $(X={\rm Spec}(R),\omega)$. Then we can define the moment map 
\[
\mu\colon  X\to \mathfrak{t}^*_{\mathbb C} \ \  (\mathfrak{t}={\rm Lie}(\mathbb T)) \, .
\] 
We consider the categorical quotient $Y=\mu^{-1}(0)/\!\!/ \mathbb{T}$. It is known that over the smooth locus $Y^{\rm sm}$,  $\omega$ descends to a symplectic form $\omega_{Y^{\rm sm}}$   (see \cite{Hitchin-hamilton}). We make the following conjecture.
\begin{problem}\label{conj-descent}
Assume the normalized minimizer of $x\in X$ is given by a conical $\mathbb {C}^*$-action with weight $\ell>0$ on $\omega$ commuting with $\mathbb T$, and the extension $(Y,\omega_Y)$ of $(Y^{\rm sm}, \omega_{Y^{\rm sm}})$ is a symplectic singularity. Show that the minimizer on $Y$ is given by the conical $\mathbb {C}^*$-action descended from the conical $\mathbb {C}^*$-action on $X$.
\end{problem} 
This conjecture is inspired by the analogue of descending a conical hyperk\"ahler metric from $X$ to $Y$. The assumption that $(Y^{\rm sm},\omega_{Y^{\rm sm}})$ extends to a symplectic singularity holds in many natural cases (see e.g. \cites{HSS-reduction, BS-quiver}). 
\end{say}
In fact, Problem \ref{conj-descent} suggests the following.

\begin{problem} \label{p:wt symmetry}
Let $(X={\rm Spec}(R),\omega)$ be a conical symplectic singularity with the conical action induced by the volume minimizer $\ord_E$. Let $\mathbb{T}$ be a maximal torus in the group of graded Hamiltonian Poisson automorphisms. Show that for any $\chi \in {\rm Hom}(\mathbb T,\mathbb G_m)$, 
\begin{equation}\label{eq-duality}
\dim R_{d,\chi}= \dim R_{d,-\chi} \, ,
\end{equation}
where $R_{d,\chi}=\{f\in R_d \mid t\cdot f= \chi(t)f, \  t\in \mathbb{T} \} $.
\end{problem}

To see the root of \eqref{eq-duality}, when the ground field is $\mathbb{C}$,  the following should be true.
\begin{problem}
Prove that there is a quaternionic structure on $R$ compatible with $\mathbb T$ (see \cite[3.7]{Etingof-quaternion}). More precisely, show that there is a conjugation $\rho$, i.e. a $\mathbb C$-anti-linear degree-preserving Poisson automorphism $\rho\colon  R \to  R$ which commutes with the maximal compact real torus $\mathbb T_c\subset \mathbb T(\mathbb C)$, such that 
\[
\rho^2(f)=\exp(\frac{2d\pi i}{r})\cdot f \mbox{\ \  for $f\in R_d$}\, ,
\]
where $r$ is the weight of the symplectic form.
\end{problem}
Assuming the weight assertion in Problem \ref{prob-contactandsym}(2) holds, $\exp(\frac{2d\pi i}{r})$ is $-1$ or $1$.
Using the real structure, one can see there is an antilinear isomorphism $\rho\colon R_{d,\chi}\to R_{d,-\chi}$.  In fact,  let \(f\in R_{d,\chi}\). By our assumption, \(\rho\) commutes with the $\mathbb T_c$-action and is
\(\mathbb C\)-antilinear, we have
\[
t\cdot \rho(f) = \rho(t\cdot f) = \rho(\chi(t)f) = \overline{\chi(t)}\,\rho(f) \, .
\]
But \(\chi(t)\in S^1\), because \(t\in \mathbb T_c \). Hence
$\overline{\chi(t)}=\chi(t)^{-1}$. Therefore
$t\cdot \rho(f)=\chi(t)^{-1}\rho(f)$,  which means that
$\rho(f)\in R_{d,-\chi}$.

There is an even more optimistic view of  proving \eqref{eq-duality}. 
The existence of a hyperk\"ahler cone structure would imply that the connected component of the graded Poisson automorphism group $G_R={\rm Aut}_{\rm Poiss}^{\mathbb G_m}(R)^{\circ}$ is reductive (see \cite[Page 15]{Etingof-quaternion}). 
 It is  known that ${\rm Lie}(G_R)=R_{r}$, whose Lie algebra structure is induced by the Poisson structure on $R$. Then for every $d$, $R_d$ is an $R_r$-module.

\begin{problem}Under the assumptions of Problem \ref{p:wt symmetry}, prove the following statements. 
\begin{enumerate} 
\item 
Show that $R_{r}$ is reductive. Therefore, it has a Chevalley involution $C\colon R_r\to R_r$ with respect to $\mathfrak{t}={\rm Lie}(\mathbb {T})\subset {R}_{r}$, i.e. $C|_{\mathfrak{t}}\colon t\to -t$.
\item (Weak) Show that for every $d$, there is a $R_r$-module isomorphism $\tau_d\colon R_d\cong R_d^C$, where $R_d^C$ is the $C$-twist of $R_d$ as a $R_r$-module.
\item (Strong) Show that there is a graded Chevalley involution $\tau\colon R\to R$, which is compatible with $C$.
\end{enumerate}
\end{problem}
If such $\tau_d$ exists, clearly $\tau_d$ induces an isomorphism $R_{d,\chi}\cong R_{d,-\chi} $. 

Another hyperk\"ahler package we can conjecture, following \cites{BPR-shortproduct, Etingof-quaternion}, is the existence of a short star-product (see {\it ibid}. for related definitions). 
\begin{problem}\label{BPR-shortproduct}
Let $(X={\rm Spec}(R),\omega)$ be a conical symplectic singularity. If we consider the degree given by the normalized volume minimizer $\frac{2}{r}\ord_E$. Show $R  $ admits a short star-product. 
\end{problem}

\subsection{Local volumes}
It is proved in \cite{XZ-localbound} that K-semistable cones whose volume is uniformly bounded from below by a positive number, are bounded. In particular, it implies that
\begin{theorem}[\cite{XZ-localbound}]\label{thm-localbounded}
Fix a finite set $I$ and a positive integer $n$. Then the set of all normalized volumes
\[
{\rm Vol}_{n,I}=\{\widehat{\rm vol}(x,X,\Delta)\mid \dim(X)=n, \, {\rm Coeff}(\Delta)\subset I\}
\]
is discrete away from 0. 
\end{theorem}
It is also known that the maximal value of ${\rm Vol}_{n,I}$ is $n^n$, realized only by the smooth point $x\in X$ ($\Delta=0$). In dimension 3, using the classification of singularities, we know more: in particular, when $\Delta=0$, the second largest volume is given by the $A_1$-singularity. See \cite{LX-cubic} and also \cite{Liu-threefolds} for further results. In general dimension, we can ask

\begin{problem}\label{p-ODP}
 If $x\in X$ is singular, then $\widehat{\rm vol}(x, X)\le 2(n-1)^n$ and the equality holds only for $A_1$-singularity.  
 \end{problem}
 An affirmative answer to this question would imply that the K-moduli of cubic hypersurfaces is the same as the GIT moduli (see \cite{SS-degree4,LX-cubic}). More generally, a better understanding of local volumes of singularities will be useful for explicit description of K-moduli spaces by the continuity method using the local-to-global inequality \cite{Liu-singularKE}. For additional works along this line, see \cite{OSS-surface, AdVL-wall1,AdVL-wall2,AdVL-wall3, Liu-threefolds} etc..

By induction on dimension and Theorem \ref{thm-SDC}, it suffices to treat isolated K-semistable log Fano cone singularities $o\in C$. Recently, there has been an important conceptual progress on this question, namely, in \cite{LM-second} the global analogue of this question has been answered affirmatively: the K-semistable Fano manifolds with the second largest volume are the smooth quadric hypersurface $X$ and the product $\mathbb{P}^{n-1}\times \mathbb{P}^1$; in other words, the only K-semistable Fano manifold with volume larger than $2n^n$ is $\mathbb{P}^n$ with $(-K_{\mathbb P^n})^n=(n+1)^n$. This in particular implies Problem \ref{p-ODP} holds for cone over smooth K-semistable Fano manifolds (see \cite[Theorem 4.15]{LM-second}).

\section{Miscellaneous Questions}

\subsection{MMP and local volume}

Recently, there has been progress of using K-stability theory to study minimal model program theory in \cite{HQZ-boundednessMMP}, where normalized volume is used as a tool to control the singularities appearing in the minimal model program sequence. 

\begin{say}[Log general type MMP]
Let $V_{\bullet}$ be a graded linear system containing an ample series on a projective klt pair $(X,\Delta)$. One can define 
\[
\alpha(X,\Delta;V_{\bullet})=\lim_{m\to \infty} m\cdot \inf\{{\rm lct}(X,\Delta; D)\,|\,D\in |V_m|\}  \, .
\]
Then there is a local-to-global inequality: for any $x\in X$, it is proved in \cite{HQZ-boundednessMMP}
\begin{equation}\label{ineq-localglobal}
\widehat{\rm vol}(x, X,\Delta) \ge  \alpha(X,\Delta;V_{\bullet})^n \cdot  {\rm vol}(V_{\bullet})\, .
\end{equation}

If $(X,\Delta)$ is a projective klt pair of log general type, i.e. $K_X+\Delta$ is big. Denote by $v={\rm vol}(K_X+\Delta)>0$ and $\alpha= \alpha(X,\Delta;K_X+\Delta)>0$. Let $X\dasharrow X'$ be a $(K_X+\Delta)$-nonnegative birational contraction. 
Then
\[
{\rm vol}(K_{X'}+\Delta')=v \mbox{ \ and \ }\alpha(X',\Delta';K_{X'}+\Delta') \ge \alpha \, .
\]
Therefore, for any $x'\in (X',\Delta')$, \eqref{ineq-localglobal} implies
\[
\widehat{\rm vol}(x', X',\Delta')\ge  \alpha^n\cdot v\,,
\] in other words, the local volume $\widehat{\rm vol}(x', X',\Delta')$ has a uniform positive lower bound. 

One can draw many conclusions from the lower bound of  the normalized volume of a klt singularity as in \cites{XZ-uniqueness, XZ-localbound}, such as an upper bound on the Cartier index, or more generally, the order of the local fundamental group; the upper bound of the log discrepancy etc..

In particular, a famous result of Shokurov says termination of any sequence of of flips follows from ACC and low-semicontinuity  of minimal log discrepancies (mld). When $(X,\Delta)$ is log general type, one can replace ACC of mld by the above argument.
\end{say}

So it is natural to ask the following question as an intermediate step to attack termination of flips.

\begin{problem}\label{prob-boundedvolume} 
Fix a klt pair $(X,\Delta)$. Is there a positive number $\varepsilon$, such that for any $(K_X+\Delta)$-nonnegative birational map $X\dasharrow X'$,  and any $x\in X'$, 
\[
\widehat{\rm vol}(x, X',\Delta') \ge \varepsilon \, . 
\]
\end{problem}

\begin{rem}
It might be too optimistic to expect Problem \ref{prob-boundedvolume} to have a positive answer when $K_X+\Delta$ is not pseudo-effective, but we don't have an example to disprove it.
\end{rem}

\subsection{Higher rank finite generation}
A novel part of K-stability theory, beyond the traditional higher dimensional geometry, is to consider valuations which are more general than divisorial valuations. A prominent class is given by quasi-monomial valuations. It also provides new tools to study problems for graded ideal sequences. 

\subsubsection{Higher rank valuations}

Let $f\colon (Y,\Delta)\to X$ be a Fano type morphism over a germ $x\in X$, and $D$ a $\mathbb{Q}$-complement. Let ${\rm LCP}(Y,\Delta+D)$ be the set of all valuations which are lc places, i.e. valuations $v\in \Val_{X,x}$ such that $A_{Y,\Delta+D}(v)=0$, which is a cone in ${\rm Val}_{X,x}$ over the dual complex $\mathcal{DR}(Y,\Delta+D)$. 
We have the following two facts concerning the divisorial valuations in ${\rm LCP}(Y,\Delta+D)$:
\begin{enumerate}
\item By \cite{BCHM}, $A:= \bigoplus_{m\in\mathbb{N}} f_*\mathcal{O}_Y(L)$ is a finitely generated $\mathcal{O}_{X,x}$-algebra for any Cartier divisor; moreover, for any divisorial valuation $v=\ord_E$   in ${\rm LCP}(Y,\Delta+D)$, the graded ring ${\rm Gr}_v(A)$ is finitely generated.
\item  When $f={\rm id}$, there exists a Koll\'ar component $E$ in ${\rm LCP}(Y,\Delta+D)$ (see \cite{Xu-kollar}).
 \end{enumerate}
 The study of K-stability inspired the more challenging question to study finite generation for higher rank quasi-monomial valuations.

\begin{say}\label{say-twotype}(\cite[Section 6]{LXZ-HRFG})
Let $Y=\mathbb{P}^2$ and $f=zx^2+zy^2+y^3$, so $D:={\rm div}(f)$ is the nodal cubic curve.   Denote by $R=\bigoplus_{m\in\mathbb{N}} H^0(Y,-mK_Y)$. 
Then
\[
{\rm LCP}(Y,D)=\left\{v_{s,t}\,|\, (s,t)\in \mathbb{R}^2_{\ge 0}\setminus(0,0) \right\} \, ,
\]
where $v_{s,t}$ is a quasi-monomial valuation with weights $(s, t)$ in the blowup of $[0 : 0 : 1]\in \mathbb{P}^2$ with respect to the two branches of $D$, with the understanding that and $v_{\lambda,0} = v_{0,\lambda} = \lambda \cdot \ord_D$. 


Then for $v_{s,t}\in {\rm LCP}(Y,D)$, ${\rm Gr}_vR$  is finitely generated if and only if 
\begin{enumerate}
\item[Type I]: ${\rm rank}_{\mathbb{Q}}(v)=1$, equivalently $\frac{t}{s} \in \mathbb{Q}$, or
\item[Type II]:  ${\rm rank}_{\mathbb{Q}}(v)=2$ with
\[
 \left\{v_{s,t} \, \Big| \,  \frac{t}{s}\in \left(\frac{7-3\sqrt{5}}{2}, \frac{7+3\sqrt{5}}{2}\right) \setminus \mathbb{Q} \right\} \, .
\]
\end{enumerate}
Similar situation also happens for the cones over other smooth del Pezzo surfaces $Y$, see \cite{Peng-degeneration}. 
\end{say}

The  higher rank quasi-monomial valuations of Type II in the above dichotomy share a stronger property as they are all special in the following sense.
\begin{say}[Special valuation]
Let $f\colon Y\to X={\rm Spec}(R)$ be a projective morphism, where $Y$ is normal, and $\Delta$ is an effective $\mathbb Q$-divisor on $Y$. We say a valuation $v$ is \emph{special} over $(Y,\Delta)$, if there is a dlt pair $(Z,\Delta_Z)$ with a birational proper morphism $g\colon Z\to Y$, such that $\Delta_Z\ge g_*^{-1}\Delta$, $-K_Z-\Delta_Z$ is ample over $X$ and $v$ is a quasi-monomial valuation in ${\rm QM}(U,E|_U)$, where $E=\lfloor \Delta_Z \rfloor$ and $U$ is an open set that intersects every lc center of $(Z,\Delta_Z)$ such that $(U,E|_U)$ is simple normal crossing. (Since ${\rm QM}(U,E|_U)$ does not depend on the choice of $U$, we denote it by ${\rm QM}(Z,\Delta_Z)$.)

The above setting can specialize to two situations which we are mostly interested in: when $X$ is a point and $(Y,\Delta)$ is a log Fano pair or $f={\rm id}$ and $(Y,\Delta)$ is a klt pair. In the latter case, a special valuation is also called a \emph{Koll\'ar valuation}.
\end{say}

Various minimizing valuations that appear in K-stability theory, e.g. valuations computing $\delta(X,\Delta)$ (assuming $\delta(X,\Delta)<\frac{\dim X+1}{\dim X}$), or minimizers of $\widehat{\rm vol}$ in ${\rm Val}_{X,x}$, are shown to be special (see \cite{LXZ-HRFG, XZ-HRFG}).

\begin{theorem}[\cite{LXZ-HRFG,XZ-HRFG,Chen-HRFG}]\label{thm-HRFG}
Any special valuation $v$ satisfies that ${\Gr}_v A$ is finitely generated, and ${\rm Spec}({\rm Gr}_v A)$ is klt. 
\end{theorem}

\begin{say}
Fix an (ordered) collection of divisorial valuations $D_1,\dots,D_k$. For $\overline{\lambda}=(\lambda_1,\dots,\lambda_k)$, one can define 
\[
\cF^{\overline{\lambda}}R= \sum_{\overline{\lambda}'\ge \overline{\lambda}} \cF^{\lambda'_1}_{D_1}\cap \cdots\cap \cF^{\lambda'_k}_{D_k}R \, ,
\]
where $\overline{\lambda}'=(\lambda'_1,\dots,\lambda_k')\ge \overline{\lambda}=(\lambda_1,\dots,\lambda_k)$ is induced by the lexicographic order. It is then clear 
this is a multiplicative decreasing filtration. As a result, we can define the  (multi-graded) ring
\[
{\rm Gr}_{D_1,\dots, D_k} R:= \bigoplus_{ \overline{\lambda}} \cF^{ \overline{\lambda} } R/  \cF^{>\overline{\lambda}}R \, .
\] 
One can check if $k=2$, ${\rm Gr}_{D_1, D_2} R\cong {\rm Gr}_{D_2,D_1} R$. In fact, we have
\begin{eqnarray*}
{\rm Gr}^{\lambda_1,\lambda_2}_{D_1, D_2} R \cong  &(\cF_{D_1}^{\lambda_1}R\cap \cF_{D_2}^{\lambda_2}R)+\cF_{D_1}^{>\lambda_1}R/ (\cF_{D_1}^{\lambda_1}R\cap \cF_{D_2}^{>\lambda_2}R)+\cF_{D_1}^{>\lambda_1}R
 \\\cong &  (\cF_{D_1}^{\lambda_1} R\cap \cF_{D_2}^{\lambda_2} R) / (\cF_{D_1}^{\lambda_1} R\cap \cF_{D_2}^{\lambda_2}R) \cap ((\cF_{D_1}^{\lambda_1}R\cap \cF_{D_2}^{>\lambda_2}R)+\cF_{D_1}^{>\lambda_1}R) \\
\cong &  (\cF_{D_1}^{\lambda_1}\cap \cF_{D_2}^{\lambda_2}) / (\cF_{D_1}^{\lambda_1}\cap \cF_{D_2}^{>\lambda_2})+ ((\cF_{D_1}^{\lambda_1}\cap \cF_{D_2}^{\lambda_2})\cap \cF_{D_1}^{>\lambda_1}) \\
\cong & (\cF_{D_1}^{\lambda_1}\cap \cF_{D_2}^{\lambda_2}) / (\cF_{D_1}^{\lambda_1}\cap \cF_{D_2}^{>\lambda_2})+ (\cF_{D_1}^{>\lambda_1}\cap \cF_{D_2}^{\lambda_2}) \\
\cong  &{\rm Gr}^{\lambda_2,\lambda_1}_{D_2, D_1} R  .
\end{eqnarray*}
Nevertheless, when $k>2$, this definition of ${\rm Gr}_{D_1,\dots, D_k} R$ in general depends on the order.

Somewhat surprisingly, in \cite{Xu-multideg}, it is proved that for any collection of divisorial valuations $D_1,\dots,D_k\in \mathcal{DR}(Y,\Delta+D)$, 
${\rm Gr}_{D_1,\dots, D_k} R $ is well defined and only depends on the set of components (but not the order of $D_i$). 

Then the key question is  when the degeneration  ${\rm Gr}_{D_1,\dots, D_k} R$ is integral, since if so the graded ring will arise as an associated graded ring of some valuation, which can be identified with a corresponding quasi-monomial valuation in ${\rm QM}(Y,\Delta)$.

This is a subtle question. It is proved in \cite{XZ-HRFG} that if $f\colon Y\to X$ is a birational morphism, $(Y,\Delta_Y)$ is a dlt pair and $-(K_Y+\Delta_Y)$ is $f$-ample, then ${\rm Gr}_{D_1,\dots, D_k} R$ is integral, where $D_1,\dots,D_k$ are the components of $\lfloor \Delta_Y \rfloor$. The conditions are somewhat relaxed in \cite{Chen-HRFG}, which also contains a more elegant proof of the integrality of ${\rm Gr}_{D_1,\dots, D_k} R$.
\end{say}

We can ask the following.
\begin{problem}
Let $D$ be a complement of a Fano type morphism $f\colon (Y,\Delta)\to X$. What is the structure of all special valuations ${\rm LCP}^{\rm S}(Y,\Delta+D)$ in ${\rm LCP}(Y,\Delta+D)$?
\end{problem}
Generally, little is known. In  \cite{LX-kollarnote}, the connectedness of ${\rm LCP}^{\rm S}(Y,\Delta+D)$ is proved when $Y=X$, using the uniqueness of minimizer of normalized volume functions.  The question  has connection with other fields, e.g. mirror symmetry.

\subsubsection{Graded ideal sequence}

\begin{say}
It is more challenging to study valuations computing the log canonical threshold of a graded sequence of ideals $\mathfrak{a}_{\bullet}$. Let $x\in (X=\Spec(R),\Delta)$ be a germ of a klt singularity. We assume $\{\fa_{\bullet}\}$ consists of $\fm_x$-primary ideals and $c={\rm lct}(X,\Delta;\fa_{\bullet})<+\infty$.  Let
\[
{\rm LCP}(X,\Delta+c\cdot \fa_{\bullet}):=\left\{v\in {\rm Val}_{X,x}\, \Big|\,  c=\frac{A_{X,\Delta}(v)}{v(\fa_{\bullet})}\right\} \, .
\]
 It is first shown in \cite{JM-quasimonomial} that  ${\rm LCP}(X,\Delta+c\cdot \fa_{\bullet})\neq \emptyset$.  Applying  \cite{Bir-BABI} to investigate an approximating sequence  $\{v_k\}_{k=1,2,\dots}$ where $v_k$ computes ${\rm lct}(X,\Delta;\frac{1}{k}\fa_k)$, it is proved in \cite{Xu-quasimonomial} that ${\rm LCP}(X,\Delta+c\cdot \fa_{\bullet})$ contains a quasi-monomial valuation, which confirms a conjecture of \cite{JM-quasimonomial}. However, many more structural questions remain to be answered.
\end{say}

\begin{problem}In the above notation, prove the following.
\begin{enumerate}
\item (Weaker version \cite{JM-quasimonomial}) All valuations in ${\rm LCP}(X,\Delta+c\cdot \fa_{\bullet})$ are quasi-monomial.
\item (Stronger version) There exists a log smooth model $(Y,E)\to X$ such that 
\[
{\rm LCP}(X,\Delta+c\cdot \fa_{\bullet}) \subseteq {\rm QM}(Y,E)
\] 
is a cone over a connected convex set.
\end{enumerate}
\end{problem}

In \cite{JM-quasimonomial}, similar questions are asked for higher jumping numbers of $\fa_{\bullet}$, but we lack the tools to attack the problem. We can also ask whether valuations with stronger properties exist in ${\rm LCP}(X,\Delta+c\cdot \fa_{\bullet})$ 
\begin{problem}In the above notation, prove the following.
\begin{enumerate}
\item (Weaker version) ${\rm LCP}(X,\Delta+c\cdot \fa_{\bullet})$ contains a valuation $v$ such that the graded algebra ${\rm Gr}_vR$ is finitely generated. 
\item (Stronger version) ${\rm LCP}(X,\Delta+c\cdot \fa_{\bullet})$ contains a Koll\'ar valuation.
\end{enumerate}
\end{problem}

\subsection{Positive characteristics} 
We could define invariants of a similar flavor in positive characteristics. However,  unlike the normalized volume, a major difference is that a valuative viewpoint (e.g. in \eqref{df-volume}) is missing in positive characteristics. 
\begin{say}
Let $\mathfrak{a}$ be an $\mathfrak{m}_x$-primary ideal on a strongly $F$-regular singularity $x\in (X={\rm Spec}(R),\Delta)$.  One can define an analogue of the log canonical threshold (known as the $F$-pure threshold) as
\[
{\rm fpt}(X,\Delta;\mathfrak{a})=\sup\{t\ge 0\mid (X,\Delta+t\cdot\mathfrak{a}) \mbox{ is strongly $F$-regular}\} \, .
\]
Then one can define the  \emph{Frobenius volume} (see \cite{LLX-localsurvey, LP-Fvolume})
\[
{\rm Fvol}(x, X,\Delta)=\inf_{\mathfrak{a}} \mult(\mathfrak{a})\cdot {\rm fpt}(X,\Delta;\mathfrak{a})^n\,.
\]
and 
the  \emph{Hilbert-Kunz Frobenius volume}
\[
{\rm Fvol_{HN}}(x, X,\Delta)=\inf_{\mathfrak{a}} \mult_{\rm HN}(\mathfrak{a})\cdot {\rm fpt}(X,\Delta;\mathfrak{a})^n\,,
\]
where in the infima run through all $\mathfrak{m}_x$-primary ideals and $ \mult_{\rm HN}(\mathfrak{a})$ is the Hilbert-Kunz multiplicity which is defined as
\[
 \mult_{\rm HN}(\mathfrak{a})=\lim_{e\to+\infty}\frac{{\rm length}(R/a^{[p^e]})}{p^{e\dim(X)}} \mbox{\ \ and \ \ } a^{[p^e]}=(f^{p^e} \mid f\in \mathfrak{a}) \, .
 \]
The comparison with the characteristic 0 definition is via an alternative interpretation of volume by \cite{Liu-singularKE}:
\[
\widehat{\rm vol}(x,X,\Delta)=\inf_{\mathfrak{a}} \mult(\mathfrak{a})\cdot {\rm lct}(X,\Delta;\mathfrak{a})^n \, .
\]
\end{say}

Let $x\in(X, \Delta)$ be a klt singularity over characteristic 0. For any ideal $\mathfrak{a}$ on $X$, let $x_p \in (X_p,\Delta_p)$ and $\mathfrak{a}_p$ be its reduction mod $p> 0$. We have \cites{HY-test=multiplier}
\[
{\rm lct}_x(X,\Delta;\mathfrak{a})=\lim_{p\to\infty} {\rm fpt}_{x_p}(X_p,\Delta_p;\mathfrak{a}_p) \, .
\]
So it is natural to ask the following
\begin{problem}[{\cite[6.19]{LLX-localsurvey}}]
Do we always have
\begin{equation}\label{eq-volumeconverge}
\widehat{\rm vol}(x,X,\Delta)=\lim_{p\to\infty}{\rm Fvol}(x_p,X_p,\Delta_p) \, ?
\end{equation}
\end{problem}
In view of \ref{say-Falpha}, it might be too optimistic to expect \eqref{eq-volumeconverge} holds. On the other hand, similarly to the conjecture in \cite{CRST-finitedegree}, we should at least expect that $\liminf_{p\to +\infty}{\rm Fvol}(x_p,X_p,\Delta_p)$ has a positive lower bound. We can ask whether there exist constants $c_n$, $C_n$ which only depend on $n=\dim(X)$ such that 
\begin{eqnarray}
c_n\cdot \widehat{\rm vol}(x,X,\Delta)\le \liminf_{p\to\infty}{\rm Fvol}(x_p,X_p,\Delta_p)\le \limsup_{p\to\infty}{\rm Fvol}(x_p,X_p,\Delta_p) \le C_n\cdot \widehat{\rm vol}(x,X,\Delta)\,.
\end{eqnarray}

\begin{say}
A better studied invariant for a strongly $F$-regular singularity is its \emph{$F$-signature} $s(X)$. It is proved in \cite{LP-Fvolume} that $s(X)$ is bounded by ${\rm Fvol}(x,X)$ from both sides
up to constants depending only on the dimension. 
\end{say}
\begin{say}[Finite degree formula]
In general, a `volume' like function $V$ of a singularity $x\in X$ should be considered to measure the volume of the `link' of $x\in X$. Therefore, it should satisfy a \emph{finite degree formula} for any finite morphism which is \'etale outside $x$. We can often expect a stronger version, namely, if $f\colon (y\in Y)\to (x\in X)$ is a finite morphism which is \'etale in codimension one, i.e  $f^*(K_X)=K_Y$, then it should satisfy
\begin{eqnarray}\label{eq-finitedegree}
V(y,Y)={\rm deg}(f)\cdot  V(x,X)\, .
\end{eqnarray}
This is established for the normalized volume \cite{XZ-uniqueness} and $F$-signature \cite{CRST-finitedegree}.  (We note that for the normalized volume, in the log case, if $f\colon (y,Y,\Delta_Y)\to (x,X,\Delta)$ with $f^*(K_X+\Delta)=K_Y+\Delta_Y$, a similar equality as in \eqref{eq-finitedegree}  holds when $f$ is Galois, but has not been established in full generality.)

So if ${\rm Fvol}$ and ${\rm Fvol_{HN}}$ are the right `volume-type' notions, we should expect the following is true.
\end{say}

\begin{problem}Let $(y\in Y)\to ( x\in X)$ be a finite morphism which is \'etale in codimension one. We have
\[
{\rm Fvol}(y,Y)={\rm deg}(f)\cdot  {\rm Fvol} (x,X)  \mbox{\ \ and \ \ } {\rm Fvol_{HN}}(y,Y)={\rm deg}(f)\cdot  {\rm Fvol_{HN}} (x,X) \, .
\]
\end{problem}

One consequence of Theorem \ref{thm-SDC} is that for any klt singularity $x\in X$ any volume $\widehat{\rm vol}(x\in X)$ is an algebraic number (see \cite{XZ-HRFG}). So we can ask the following question.
\begin{problem}
Is $ {\rm Fvol} (x,X) $  an algebraic number?
\end{problem}
We do not have reason to believe the above question has an affirmative answer, but a counterexample will probably shed light on the difference between volumes and $F$-volumes.

\begin{say}\label{say-Falpha}
We can similarly define the (local) \emph{Frobenius $\alpha$-invariant} (see \cites{Pande-Falpha, ST-Fsignature})
\[
\alpha_F(x,X,\Delta)=\inf_{\mathfrak{a}} \overline{ \ord}(\mathfrak{a})\cdot {\rm fpt}(X,\Delta;\mathfrak{a}) \, ,
\]
where $ \overline{\ord}(\mathfrak{a})=\lim_{m\to \infty}\frac{\ord(\mathfrak{a}^m)}{m}$ and $\ord(\cdot)$ means order along the maximal ideal (i.e. the largest integer $r$ such that $\fm_x^r$ contains the given ideal). It is clear from the definition that $\alpha_F(x,X,\Delta)\le {\rm fpt}(X,\Delta;\mathfrak{m}_x)$, hence
\[
e(R)\cdot \alpha_F(x,X,\Delta)^n\le  {\rm Fvol} (x,X,\Delta)\, ,
\]
where $e(R)$ is the Hilbert-Samuel multiplicity of $R$. As noticed in \cite{Pande-Falpha}, in general
\[
\alpha(x,X,\Delta)\neq\lim_{p\to+\infty}\alpha_F(x_p,X_p,\Delta_p) \,.
\]
\end{say}

\bibliography{refxu} 
\bibliographystyle{ieeetr}
\end{document}